\documentclass[a4paper,10pt]{article}
\usepackage{amssymb}
\usepackage{amsmath}

\usepackage[OT2,OT1]{fontenc}
\newcommand\cyr{%
\renewcommand\rmdefault{wncyr}%
\renewcommand\sfdefault{wncyss}%
\renewcommand\encodingdefault{OT2}%
\normalfont
\selectfont}
\DeclareTextFontCommand{\textcyr}{\cyr}

\numberwithin{equation}{section}
\newcommand{\twosum}[2]{\sum_{\begin{array}{c} {\scriptstyle #1}\\
{\scriptstyle #2} \end{array}}}
\newcommand{\ep}{\varepsilon}
\newcommand{\Q}{\mathbb{Q}}
\newcommand{\N}{\mathbb{N}}
\newtheorem{theorem}{Theorem}
\newtheorem{lemma}{Lemma}
\renewcommand{\mod}[1]{\hspace{-1.5mm}\pmod{#1}}
\begin{document}
\title{The Average Analytic Rank of Elliptic Curves}
\author{D.R. Heath-Brown\\Mathematical Institute, Oxford}
\date {}
\maketitle

\begin{abstract}
All the results in this paper are conditional on the Riemann
Hypothesis for the $L$-functions of elliptic curves.  Under this
assumption, we show that the average analytic rank of all
elliptic curves over $\Q$ is at most 2, thereby improving a result of
Brumer \cite{brum}.  We also show that the average within any family
of quadratic twists is at most $3/2$, improving a result of Goldfeld
\cite{gold}. A third result concerns the density of curves with 
analytic rank at least $R$, and shows that the proportion of such 
curves decreases faster than exponentially as $R$ grows.  The proofs
depend on an analogue of Weil's ``explicit formula''.
\end{abstract}

\section{Introduction}
The purpose of this paper is to establish upper bounds for the average
of the analytic rank of elliptic curves defined over
$\mathbb{Q}$.  The article by Rubin and Silverberg \cite{RS}
gives an excellent survey of this topic.
Our first
result concerns the average over all such curves, and sharpens an estimate of 
Brumer \cite{brum}.  We introduce at the outset the minor technical trick of
counting the curves
\[E=E_{r,s}:\; y^{2}=x^{3}+rx+s\]
with a weight
\[w_{T}(E)=w_{1}(T^{-1/3}r) w_{2}(T^{-1/2}s),\]
where $w_{1},w_{2}$ are infinitely differentiable non-negative functions
of compact support, vanishing at the origin.  
We define $\Delta_{E}=-16(4r^{3}+27s^{2})$ and we
write
\[{\cal C}=\{E_{r,s}:p^{4}|r\Rightarrow p^{6}\nmid s,\; \Delta_{E}\not=0\}\] 
and
\[{\cal S}(T)=\sum_{E\in{\cal C}}w_{T}(E).\]
Our principal result is then the following.
\begin{theorem}
Assume that the $L$-functions of all the curves $E_{r,s}$ satisfy
the Riemann Hypothesis.  Then
\[\frac{1}{{\cal S}(T)}\sum_{E\in{\cal C}}w_{T}(E)r(E)\leq 2+o(1),\]
as $T\rightarrow\infty,$ where $r(E)$ is the analytic rank of $E.$
\end{theorem}
Thus the average analytic rank, taken over all elliptic curves defined
over $\Q$, is at most $2$.  This improves on the corresponding result
of Brumer \cite{brum}, in which it was shown that the average is at
most $2.3$.

We shall also investigate the proportion of elliptic curves
$E$ which have large rank.  We define two sets,
\[{\cal D}(T)=\{E_{r,s}: |r|\leq T^{1/3},|s|\leq T^{1/2}, 
\Delta_{E_{r,s}}\not=0\},\]
and
\[{\cal C}(T)=\{E_{r,s}\in{\cal D}(T): p^{4}|r\Rightarrow p^{6}\nmid
s\}.\]
We then have the following result.
\begin{theorem}
Assume that the $L$-functions of all the curves $E_{r,s}$ satisfy
the Riemann Hypothesis.  Then for any positive integer $R$ we have
\[\frac{\#\{E\in{\cal C}(T): r(E)\geq R\}}{\#{\cal C}(T)}\ll 
(3R/2)^{-R/12},\]
where the implied constant is absolute.
\end{theorem}
Thus the proportion of curves with rank $R$ decreases faster than 
exponentially.  We remark that it may be possible to improve the values of the
constants $3/2$ and $12$ which occur in the theorem.  We have merely given the 
simplest values that the method allows.

Since
\[\#{\cal C}(T)\ll T^{5/6},\]
it follows that
\[r(E)\leq 11\frac{\log T}{\log\log T}\]
for $E\in{\cal C}(T)$ and sufficiently large $T.$  Such results are 
already known (see Mestre \cite{Mestre} and
Brumer \cite{brum}).  However the fact that our theorem actually contains this
estimate demonstrates that we have achieved the best rate of decay with
respect to $R$ that one can currently hope for.

We shall also consider the set of quadratic twists
\[E_{D}: Dy^{2}=x^{3}+rx+s\] 
of a fixed elliptic curve $E$ of conductor $N$, say.  
This family has previously been 
investigated by Goldfeld \cite{gold}.  It is of some interest to
separate the odd rank twists from those of even rank.  We therefore
define $L_D(s)$ to be the $L$-function of $E_D$, and $w_D=\pm 1$ to
be the sign of the functional equation for $L_D(s)$.  Thus if $w$ is
the corresponding sign for the original curve $E$ we have
\begin{equation}\label{1.1}
w_D=w\frac{D}{|D|}\chi_D(N),
\end{equation}
for $(D,N)=1$, where $\chi_D$ is the real primitive character
associated to the quadratic field $\mathbb{Q}(\sqrt{D})$.  We then set
\[{\cal T}=\{D: (D,N)=1\},\]
and
\[{\cal T}^{\pm}=\{D\in{\cal T}: w_D=\pm 1\},\]
where $D$ is restricted to run over fundamental discriminants in each case.
It
follows that $r(E_{D})$ is even for $D\in{\cal T}^{+},$ 
and odd for $D\in{\cal T}^{-}.$

For technical reasons we find it convenient to count the twists $E_{D}$ with 
a smooth weight.  
We therefore choose a three times differentiable non-negative function
$w(x),$ supported on a compact subset of either $(-\infty,0)$ or $(0,\infty),$
and we define
\[{\cal W}^{\pm}(T)=\sum_{D\in{\cal T}^{\pm}}w(D/T).\]
Our result is then the following.
\begin{theorem}
Let $E$ be a fixed elliptic curve defined over $\Q$.  Suppose that 
the functions 
$L_D(s)$ all satisfy the Riemann Hypothesis.  Then
\[\frac{1}{{\cal W}^{\pm}(T)}\sum_{D\in{\cal T}^{\pm}}w(D/T)r(E_{D})
\leq \frac{3}{2}+o(1),\]
as $T\rightarrow\infty.$
\end{theorem}
Of course it is natural to apply this result with a weight $w$ which 
approximates to the characteristic function of an interval.
Thus within a family of quadratic twists the average analytic
rank would be at most
$\frac{3}{2}$, whether one restricts to curves of odd rank or to curves of
even rank.
This may be compared with a result of Goldfeld \cite[Proposition
2]{gold}, who considers the set
${\cal T}$ only, and in which the constant
$\frac{3}{2}$ is replaced by $\frac{13}{4}.$
The reader should note that our theorem requires $L_D(s)$ to satisfy
the Riemann Hypothesis for every integer $D$, even though the sets
${\cal T}^{\pm}$ contain only integers $D$ which are coprime to $N$.

Naturally we expect that the above results should remain true if we
replace the analytic rank $r(E)$ by the arithmetic rank, which we
denote by $R(E)$. Results of Kolyvagin \cite{Koly1}, \cite{Koly2} and
Gross and Zagier \cite{GZ} show that
\begin{equation}\label{1.2}
R(E)=\left\{\begin{array}{cc} 0, & \mbox{if }\;\; r(E)=0,\\
1,& \mbox{if }\;\; r(E)=1.\end{array}\right.
\end{equation}
Theorem 3 then has the following corollary.
\begin{theorem}
Let $E$ be an elliptic curve defined over $\Q$.  Suppose that the functions 
$L_D$ all satisfy the Riemann Hypothesis.  Then
\begin{eqnarray*}
\frac{1}{{\cal W}^{+}(T)}\hspace{2mm}\sum_{D\in{\cal T}^{+},\, 
R(E_{D})=0}w(D/T)&\geq&
\frac{1}{{\cal W}^{+}(T)}\hspace{2mm}\sum_{D\in{\cal T}^{+},\, 
r(E_{D})=0}w(D/T)\\
&\geq &\frac{1}{4}+o(1)
\end{eqnarray*}
and
\begin{eqnarray*}
\frac{1}{{\cal W}^{-}(T)}\hspace{2mm}\sum_{D\in{\cal T}^{-},\, 
R(E_{D})=1}w(D/T)&\geq&
\frac{1}{{\cal W}^{-}(T)}\hspace{2mm}\sum_{D\in{\cal T}^{-},\, 
r(E_{D})=1}w(D/T)\\
&\geq& \frac{3}{4}+o(1)
\end{eqnarray*}
as $T\rightarrow\infty.$
\end{theorem}

Thus at least $1/4$ of all curves in ${\cal T}^+$ would have
arithmetic rank 0, and at least $3/4$ of all curves in ${\cal T}^-$
would have rank 1.

{\bf Acknowledgements.}  The bulk of this work was carried out in 1991,
at the Institute for Advanced Study, Princeton.  It was
prepared for publication more recently
while the author was a guest of the American
Institute of Mathematics.  It is a pleasure to thank both these bodies
for their hospitality and financial assistance.

\section{Preliminaries}
Our starting point is the `explicit formula' in the form given by
Brumer \cite[\S 2]{brum}.  We apply this to an arbitrary elliptic
curve $E$, which is of course now known to be modular, following the
work of Wiles \cite{wiles}, Taylor and Wiles \cite{WT} and 
Breuil {\em et. al.} \cite{breuil}.  We write $N_E$ for the conductor
of $E$, and
\[L_E(s)=\sum_{n=1}^{\infty}a_n(E)n^{-s}\]
for the $L$-function of $E$.  We note that if $p\nmid N_E$ then 
\[a_{p}(E)=\alpha_{p}+\overline{\alpha}_{p}\]
where 
\[|\alpha_{p}|=\sqrt{p}.\]
For the remaining primes
$p|N_E$ the coefficients $a_p(E)$ are always $0,1$ or $-1$.  Finally
we define
\[c_{p^{k}}(E)=\begin{cases} -a_{p}(E)^{k}/kp^{k}, & p|N_{E},\\
-(\alpha_{p}^{k}+\overline{\alpha}_{p}^{k})/kp^{k}, & p\nmid N_{E}.
\end{cases}\]

We take the weight function $F(t)$ in \cite[Lemma 2.1]{brum}
to be
\[F(t)=h_X(t)=h(t/\log X)\] 
where $X\ge 2$ and
\[h(t)=\left\{\begin{array}{cc} 1-|t|, & |t|\leq 1,\\ 0, & |t|>1. \end{array}
\right.\]
We define the Fourier transform of a function $f(x)$ by
\[\hat{f}(t)=\int_{-\infty}^{\infty}e^{-2\pi ixt}f(x)dx.\]
Although this convention differs from Brumer's, the two alternative 
definitions of $\hat{f}(0)$ agree.
Since $\hat{h}(0)=h(0)=1$ and $\hat{h}(t)\geq 0$ for all real $t,$ we deduce
from the estimates of Brumer \cite[\S 2]{brum} that
\begin{equation}\label{2.1}
r(E)\leq \frac{\log N_{E}}{\log X}+
\frac{2}{\log X}(U_{1}(E,X)+U_{2}(E,X)) +O(\frac{1}{\log X}),
\end{equation}
where 
\[U_{k}(E,X)=\sum_{p^{k}\leq X,\,p\ge 5}
c_{p^{k}}(E)(\log p^{k})h_{X}(\log p^{k}).\]
In particular we have
\[U_1(E,X)=-\sum_{5\le p\leq X}\frac{\log p}{p}h_{X}(\log p)a_p(E).\]

We begin by considering $U_2(E,X)$.  For Theorem 1 we can use Brumer's
work \cite[p. 457]{brum}, which shows that
\begin{eqnarray}\label{2.2}
\sum_{E\in{\cal C}}w_{T}(E)U_2(E,X)&=&
{\cal S}(T)\{\frac{\hat{h}(0)}{4}\log X +O(T^{-1/2}X^{1/2}\log X)\nonumber\\
&&\hspace{2cm}\mbox{}+O(T^{-3/4}X^{9/10}(\log X)^{9/5})\}\nonumber\\
&=& {\cal S}(T)\{\frac{1}{4}\log X +O(1)\},
\end{eqnarray}
providing that $X\le T^{5/6-\delta}$ with a fixed $\delta>0$.

Similarly for Theorem 3 we may use the work of Goldfeld
\cite[p. 116]{gold} which produces
\begin{equation}\label{2.2a}
U_2(E_D,X)=\frac{\hat{h}(0)}{4}\log X +O(\log\log D)=
\frac{1}{4}\log X +O(\log\log D),
\end{equation}
with an implied constant depending on $E$.  Finally, for Theorem 2 we
note that
\[|c_{p^{k}}(E)|\leq 2k^{-1}p^{-k/2},\]
from the definition, whence
\begin{equation}\label{2.3}
|U_{2}(E,X)|\leq 2\sum_{5\le p\leq \sqrt{X}}\frac{\log p}{p}\leq \log X
+O(1).
\end{equation}

It remains therefore to consider the behaviour of $U_1(E,X)$, for which
we shall require slightly different techniques in each case.

Before leaving this section we need to record one further result given
by Brumer \cite[(2.13)]{brum}.  In view of Brumer's `Note added in
proof' \cite[p. 472]{brum}, we may state the result as follows.
\begin{lemma}
Let $k$ be an even $C^1$ continuous function with support in $[-1,1]$,
and suppose that $\hat{k}(t)=O_{\delta}((1+|t|)^{-1-\delta})$ for 
some $\delta>0$.  Then
if $E$ is an elliptic curve of conductor $N_E$ we have
\[\sum_{p\leq X}\frac{\log p}{p}k(\frac{\log p}{\log X})a_p(E)
\ll_{\delta}
(\log N_E)(\log X)\{||k||_{\infty}+||(1+|t|)^{1+\delta}\hat{k}||_{\infty}\}.\]
\end{lemma}
Although Brumer proves this only when $X>10\log N_E$ it is automatically
true for smaller $X$, by virtue of the bound $|a_p(E)|\le 2\sqrt{p}$.

\section{Theorem 1---Initial Transformations}

Let ${\cal D}$ be the set of all curves $E_{r,s}$, including those for
which $\Delta=0$.  For any curve $E_{r,s}\in{\cal D}$ we define
\begin{equation}\label{3.1}
\sigma_{p}(E_{r,s})=-\tau_{p}^{-1}\sum_{t,x\!\mod{p}}(\frac{t}{p})e_{p}(tx^{3}
+txr+ts),
\end{equation}
where $\tau_p$ is the usual Gauss sum.  Here we adopt the standard 
convention that $e_{p}(x)=\exp(2\pi ix/p).$  
Then, according to Brumer \cite[(3.1)]{brum}, we have
\begin{equation}\label{3.1a}
a_p(E)=\sigma_p(E)\;\;\;(p\ge 5)
\end{equation}
for every $E\in{\cal C}(T)$.
As Brumer remarks, this
formula is valid even when $E$ is an elliptic curve with
singular reduction modulo $p$.  We also observe that for every
$E_{r,s}\in{\cal D}$ we have
\[\sigma_p=-\sum_{x\!\mod{p}}(\frac{x^3+xr+s}{p}),\]
whence
\begin{equation}\label{3.1b}
|\sigma_p|\le 2\sqrt{p}.
\end{equation}

The key estimate required for Theorem 1 is then as follows.
\begin{lemma}
For any $\ep>0$ we have
\begin{eqnarray*}
\lefteqn{\sum_{P<p\leq 2P}|\sum_{E\in{\cal D}}w_{T}(E)\sigma_p(E)|}
\hspace{1cm}\\
&\ll & P^{\ep}
(P^{1/2}T^{5/6}+P^{3/2}T^{1/2}+P^{2}T^{1/6}+P^{7/2}T^{-1}).
\end{eqnarray*}
\end{lemma}

Here, and throughout the paper, we allow the constant implied by the
$\ll$ symbol to depend on $\ep$.

For convenience of notation we shall write
\[\sum_{P<p\leq 2P}|\sum_{E\in{\cal D}}w_{T}(E)\sigma_p(E)|=\Sigma.\]
We begin the proof of Lemma 2 
by observing that the value $t=0$ in (\ref{3.1}) may be
omitted, since $(\frac{0}{p})=0.$  We can then substitute $y=tx$ for
$x,$ giving
\[\sigma_p(E_{r,s})=-\tau_{p}^{-1}\sum_{t\not\equiv 0\!\mod{p}}\;\;
\sum_{y\!\mod{p}}(\frac{t}{p})e_{p}(t^{-2}y^{3}+yr+ts).\]
(Here we interpret $t^{-2}y^3$ modulo $p$.)
Hence
\[
|\sum_{E\in{\cal D}}w_{T}(E)\sigma_p(E)|\leq p^{-1/2}|\sum_{t,y}
(\frac{t}{p})e_{p}(t^{-2}y^{3})S_{1}S_{2}|,
\]
where
\[S_{1}=\sum_{r=-\infty}^{\infty}w_{1}(T^{-1/3}r)e_{p}(yr)\]
and
\[S_{2}=\sum_{s=-\infty}^{\infty}w_{2}(T^{-1/2}s)e_{p}(ts).\]
According to the Poisson summation formula the sum $S_{1},$ for
example, is
\[T^{1/3}\sum_{m=-\infty}^{\infty}\hat{w}_{1}(T^{1/3}(m+\frac{y}{p})).\]
Moreover, since $w_{1}$ has derivatives of all orders, it follows that
\begin{equation}\label{3.2}
\frac{d^n}{dx^n}\hat{w}_{1}(x)\ll_{n,A} (1+|x|)^{-A}
\end{equation}
for any real $x$, any fixed integer $n\ge 0$, and any fixed $A>0.$  
To bound the sum over $m$ it is
convenient to fix the range for $y$ so that $|y|\leq p/2.$
We then conclude that
\[\sum_{m=-\infty}^{\infty}\hat{w}_{1}(T^{1/3}(m+\frac{y}{p}))\ll
1,\]
and
\[\sum_{m\not=0}\hat{w}_{1}(T^{1/3}(m+\frac{y}{p}))\ll T^{-A/3},\]
for $|y|\leq p/2.$  The sum $S_{2}$ may be handled
similarly, and we conclude that
\[S_{1}S_{2}
=T^{5/6}\hat{w}_{1}(\frac{yT^{1/3}}{p})\hat{w}_{2}(\frac{tT^{1/2}}{p})
+O(T^{5/6-A/3}).\]
Moreover any terms for which $P/2<|y|\leq p/2$ or $P/2<|t|\leq p/2$ are 
$O(T^{5/6-A/3}).$  We therefore arrive at the estimate
\[\Sigma \ll  P^{-1/2}T^{5/6}\sum_{p}|\twosum{|t|\leq P/2}{t\not= 0}
\hat{w}_{2}(\frac{tT^{1/2}}{p})(\frac{t}{p})\Sigma_{1}(t,p)|
+P^{5/2}T^{5/6-A/3},\]
where
\[ \Sigma_{1}(t,p)=\sum_{|y|\leq P/2}\hat{w}_{1}(\frac{yT^{1/3}}{p})
e_{p}(t^{-2}y^{3}).\]
If we take $A=6$ the final term is
\[P^{5/2}T^{-7/6}\le P^{7/2}T^{-1},\]
which is satisfactory for Lemma 2.

It may be worth observing at this point that the bound (\ref{3.2}), together 
with its analogue for $\hat{w}_{2},$ yields 
\[\sum_{t,y}|\hat{w}_{1}(\frac{yT^{1/3}}{p})
\hat{w}_{2}(\frac{tT^{1/2}}{p})|\ll (1+T^{-1/3}P)(1+T^{-1/2}P).\]
whence one trivially has
\[\Sigma\ll P^{1/2}(T^{1/3}+P)(T^{1/2}+P)+P^{7/2}T^{-1}.\]
This is essentially the estimate of Brumer \cite{brum}, and suffices to prove
Theorem 1 with the upper bound $2.3+o(1)$.

In order to improve on the above trivial 
argument we shall take advantage of the
averaging over $p$ to show that the oscillating term
$e_{p}(t^{-2}y^{3})$ provides some cancellation.  In order to do this we
first replace $p$ by a new variable $k$ which runs over all integers,
both prime and composite, weighted by a function $w_{3}(k/P),$ where
$w_{3}$ is an infinitely differentiable non-negative function, supported
on $[\frac{1}{2},\frac{5}{2}]$ and strictly positive on $[1,2].$ 
We then have
\begin{eqnarray*}
\lefteqn{\sum_{p}|\hspace{-2mm}\twosum{|t|\leq P/2}{t\not= 0}
\hat{w}_{2}(\frac{tT^{1/2}}{p})(\frac{t}{p})\Sigma_{1}(t,p)|}
\hspace {2cm}\\
& \ll &
\sum_{k}w_{3}(kP^{-1})|\hspace{-2mm}\twosum{0<|t|\leq P/2}{(t,k)=1}
\hat{w}_{2}(\frac{tT^{1/2}}{k})(\frac{t}{k})\Sigma_{1}(t,k)|.
\end{eqnarray*}
Here we define the Jacobi symbol $(t/k)$ to be zero whenever $k$ is
even. We now wish to bring the summation over $t$ outside the modulus signs. 
In order to do this we observe that
\[\max_{P/2\leq k\leq 5P/2}|\hat{w}_{2}(\frac{xP}{k})|
\ll_{A} (1+|x|)^{-A}\]
for any $A>0.$  Hence, on defining
\begin{equation}\label{3.3}
M(t)=\min\{1\, ,\, (|t|T^{1/2}P^{-1})^{-1}\},
\end{equation}
we see that
\[\max_{P/2\leq k\leq 5P/2}|\hat{w}_{2}(\frac{tT^{1/2}}{k})|\ll_A 
M(t)^{A},\]
for any $A>0.$  We now have
\begin{eqnarray*}
\lefteqn{\sum_{k}w_{3}(kP^{-1})|\hspace{-2mm}\twosum{0<|t|\leq P/2}{(t,k)=1}
\hat{w}_{2}(\frac{tT^{1/2}}{k})(\frac{t}{k})\Sigma_{1}(t,k)|}\hspace{3cm}\\
& \ll & \sum_{k}w_{3}(kP^{-1})\twosum{t\not=0}{(t,k)=1}M(t)^{A}
|\Sigma_{1}(t,k)|,
\end{eqnarray*}
whence
\begin{equation}\label{3.4}
\Sigma \ll P^{-1/2}T^{5/6}\sum_{k}w_{3}(kP^{-1})
\sum_{(t,k)=1}M(t)^{A}|\Sigma_{1}(t,k)|+P^{7/2}T^{-1}.
\end{equation}

\section{Lemma 2---The Kernel of the Proof}
In order to perform the averaging over $k$ we shall use Cauchy's
inequality to reverse the order of summations in (\ref{3.4}).  In view of 
(\ref{3.3}) and the fact that $w_3(x)$ is supported on
$[\frac{1}{2},\frac{5}{2}]$, 
we have
\[\twosum{k,t}{t\not=0}w_{3}(kP^{-1})|t|^{-1}M(t)^{A}
\ll P\log P.\]
On applying
Cauchy's inequality to (\ref{3.4}) we deduce that
\begin{equation}\label{4.1}
\Sigma \ll P^{-1/2}T^{5/6}\{P\log P\}^{1/2}\Sigma_{1}^{1/2}
+P^{7/2}T^{-1},
\end{equation} 
where
\begin{eqnarray*}
\Sigma_{1}&=&\sum_{k}w_{3}(kP^{-1})
\sum_{(t,k)=1}|t|M(t)^{A}.\,
|\Sigma_{1}(t,k)|^{2}\\
&=&\sum_{t}|t|M(t)^{A}
\sum_{|y_{1}|,|y_{2}|\leq P/2}\hspace{2mm}\sum_{(k,t)=1}w_{4}(kP^{-1})
e_{k}(t^{-2}\{y_{1}^{3}-y_{2}^{3}\}),
\end{eqnarray*}
with
\[w_{4}(x)=w_{3}(x)\hat{w}_{1}(\frac{y_{1}T^{1/3}}{Px})
\overline{\hat{w}_{1}(\frac{y_{2}T^{1/3}}{Px})}.\]
Thus $w_{4}$ is supported in $[\frac{1}{2},\frac{5}{2}],$ and for any 
fixed integer $n\ge 0$ and any $A>0$ we have
\begin{equation}\label{4.2}
\frac{d^{n}}{dx^{n}}w_{4}(x)\ll_{A,n} m(y_{1})^{A}m(y_{2})^{A},
\end{equation}
by (\ref{3.2}), where
\begin{equation}\label{4.3}
m(y)=\min\{1\, ,\, (\frac{|y|T^{1/3}}{P})^{-1}\}.
\end{equation}

We now write 
\[a=y_{1}^{3}-y_{2}^{3},\; b=t^{2}.\]
If $b\overline{b}\equiv 1\mod{k}$ and $k\overline{k}\equiv 1\mod{b},$
then $b\overline{b}+k\overline{k}\equiv 1\mod{bk},$ for $(b,k)=1$, so that
\[\frac{\overline{b}}{k}+\frac{\overline{k}}{b}-\frac{1}{bk}\]
is an integer.  Thus $e_{k}(a\overline{b})=e_{bk}(a)e_{b}(-a\overline{k}).$
With this in mind we define
\[\rho(x)=w_{4}(x)e(\frac{a}{bPx}),\]
so that
\[\Sigma_{1} =
\sum_{t,y_{1},y_{2}}|t|M(t)^{A}
\sum_{(k,t)=1}\rho(kP^{-1})e_{b}(-a\overline{k}).\]
We decompose the inner sum into residue classes modulo $b,$ and apply the
Poisson summation formula to obtain
\begin{eqnarray}\label{4.4}
\twosum{j\!\mod{b}}{(j,b)=1}e_{b}(-a\overline{j})
\sum_{k\equiv j\!\mod{b}}\rho(kP^{-1}) &=&\sum_{j}e_{b}(-a\overline{j})
\sum_{m=-\infty}^{\infty}\rho(\frac{j+bm}{P})\nonumber\\
&\hspace{-4cm}=&\hspace{-2cm}\sum_{j}e_{b}(-a\overline{j})
\sum_{n=-\infty}^{\infty}e_{b}(nj)\frac{P}{b}\hat{\rho}(\frac{nP}{b}).
\end{eqnarray}
At this point we observe that
\[\hat{\rho}(x)\ll x^{-2}\sup_v |\rho''(v)|\]
and that
\[\rho''(x)\ll (1+|\frac{a}{bP}|^{2})m(y_{1})^{A}m(y_{2})^{A}\]
by (\ref{4.2}),
whence
\[\hat{\rho}(\frac{nP}{b})\ll
n^{-2}(\frac{b^{2}}{P^{2}}+\frac{a^{2}}{P^{4}})m(y_{1})^{A}m(y_{2})^{A}.\]
However
\begin{eqnarray*}
\frac{b^{2}}{P^{2}}+\frac{a^{2}}{P^{4}}&\ll &
P^{2}T^{-2}\{(\frac{|t|T^{1/2}}{P})^{4}+(\frac{|y_{1}|T^{1/3}}{P})^{6}+
(\frac{|y_{2}|T^{1/3}}{P})^{6}\}\\
&\ll& P^{2}T^{-2}M(t)^{-4}m(y_{1})^{-6}
m(y_{2})^{-6},
\end{eqnarray*}
by (\ref{3.3}) and (\ref{4.3}).  It therefore follows that the terms 
$n\not=0$ in (\ref{4.4}) are
\begin{eqnarray*}
&\ll & \sum_{j\!\mod{b}}\sum_{n\not=0}\frac{P}{b}n^{-2}P^{2}T^{-2}
M(t)^{-4}m(y_{1})^{A-6}m(y_{2})^{A-6}\\
& \ll & P^{3}T^{-2}M(t)^{-4}m(y_{1})^{A-6}
m(y_{2})^{A-6}.
\end{eqnarray*}
On choosing $A=8$ we see that the contribution to $\Sigma_{1}$ is
\begin{eqnarray*}
\lefteqn{P^{3}T^{-2}\sum_{t,y_{1},y_{2}}|t|
M(t)^{4}m(y_{1})^{2}m(y_{2})^{2}}\hspace{4cm}\\
& \ll & P^{3}T^{-2}\left\{\sum_{t}|t|M(t)^{4}\right\}
\left\{\sum_{y}m(y)^{2}\right\}^{2}\\
& \ll & P^{3}T^{-2}(1+PT^{-1/2})^{2}(1+PT^{-1/3})^{2}
\end{eqnarray*}
by (\ref{3.3}) and (\ref{4.3}).  The contribution to $\Sigma$ itself is then
\begin{eqnarray*}
&\ll& P^{3/2}T^{-1/6}(\log P)^{1/2}(1+PT^{-1/2})(1+PT^{-1/3})\\
&\ll& P^{\ep}(P^{3/2}T^{-1/6}+
P^{5/2}T^{-1/2}+P^{7/2}T^{-1}),
\end{eqnarray*}
by (\ref{4.1}), and this is satisfactory for Lemma 2, since
\begin{eqnarray*}
P^{5/2}T^{-1/2}&=&\{P^{3/2}T^{1/2}.\,P^{7/2}T^{-1}\}^{1/2}\\
&\le&\max\{P^{3/2}T^{1/2}\,,\,P^{7/2}T^{-1}\}\\
&\le& P^{3/2}T^{1/2}+P^{7/2}T^{-1}.
\end{eqnarray*}

\section{Lemma 2---A Highest Common Factor Sum}
It remains to handle the terms $n=0$ in (\ref{4.4}).  Since 
\[|\hat{\rho}(0)|\leq \int_{-\infty}^{\infty}|\rho(x)|dx
= \int_{-\infty}^{\infty}|w_{4}(x)|dx\ll m(y_{1})^{A}m(y_{2})^{A},\]
by (\ref{4.2}), the contribution to $\Sigma_{1}$ is
\begin{equation}\label{5.1}
\ll \sum_{t,y_{1},y_{2}}\frac{P}{b}|t|M(t)^{A}
m(y_{1})^{A}m(y_{2})^{A}\left|
\twosum{j\!\mod{b}}{(j,b)=1}e_{b}(-a\overline{j})\right|.
\end{equation}
The sum over $j$ is a Ramanujan sum which may be evaluated as
\[\sum_{d|a,b}d\mu(b/d).\]
We therefore see that (\ref{5.1}) is 
\begin{equation}\label{5.2}
\ll P\sum_{t=1}^{\infty}t^{\ep-1}M(t)^{A}
\sum_{|y_{1}|\leq y_{2}}m(y_{2})^{A}(a,b),
\end{equation}
for any $\ep>0,$ where $(a,b)$ denotes the highest 
common factor of $a$ and $b.$  The terms with $y_1=y_2=0$ are
\[P\sum_{t=1}^{\infty}t^{\ep+1}M(t)^{A}\ll P^{3+\ep}T^{-1}+P^{1+\ep}.\]
For the remaining sum we shall use the following lemma.
\begin{lemma}
For any $U,V\geq 1$ and any $\ep>0,$ we have
\[\sum_{1\leq u\leq U}\sum_{|w|\leq v\leq V}(u^{2},v^{3}-w^{3})\ll
U^{1+\ep}V(U^{2}+V).\]
\end{lemma}

We now see that the ranges $U/2<|t|\leq U$ and $V/2<y_{2}\leq V$
contribute
\[\ll PU^{\ep-1}M(U)^{A}m(V)^{A}U^{1+\ep}V(U^{2}+V)\]
to (\ref{5.2}), and hence to $\Sigma_1$.  
We choose $A=4$, say, and sum $U$ and $V$ over powers of $2$ to
obtain a total contribution to $\Sigma_{1}$ of
\[\ll P^{3+\ep}T^{-1}+P^{1+2\ep}(1+P^{3}T^{-4/3}+P^{2}T^{-2/3}).\]
This is satisfactory for Lemma 2, by (\ref{4.1}).

It remains to prove Lemma 3.  We write $S$ for the sum to be estimated,
and we take
\[(u^{2},v^{3}-w^{3})=d=\prod p^{e}.\]
We also define
\[\delta=\prod p^{[(e+1)/2]},\]
where $[x]$ denotes the integer part of $x,$ as usual.  It follows that
$\delta|u$ whenever $d|u^{2},$ so that $u$ takes at most $U/\delta$ 
values for each given value of $d$.  We therefore see that
\[S\leq U\sum_{d\leq U^{2}}d\delta^{-1}\#\{v,w: d|v^{3}-w^{3}\}.\]
We now consider the value of $(d,v^{3})$ which we denote by 
$\alpha=\prod p^{f}.$  On defining
\[\beta=\prod p^{[(f+2)/3]},\]
we see that $\beta|v,$ and since $(d,v^{3})=(d,w^{3})$ we also must
have $\beta|w.$  We may therefore write $v=\beta v'$ and $w=\beta w',$
whence
\begin{equation}\label{5.3}
v'^{3}\equiv w'^{3}\!\mod{\gamma},
\end{equation}
where $\gamma=\prod p^{g}$ with
\[g=\max\{e-3[\frac{f+2}{3}]\, ,\, 0\}.\]
By construction we have $(v',\gamma)=1$ so that the congruence (\ref{5.3})
has at most $3^{\omega(\gamma)}\ll U^{\ep}$ solutions $w'\mod{\gamma},$
for each value of $v'.$  (Here $\omega(\gamma)$ is the number of distinct
prime factors of $\gamma.$)  Thus, for given values of $d,\alpha,\beta$ and 
$\gamma,$ there are at most $1+V/\beta$ possible choices of $v,$ to each of 
which there correspond $O(U^{\ep}(1+V/\beta\gamma))$ possible values of $w.$
We therefore conclude that
\begin{eqnarray*}
S&\ll& U^{1+\ep}\sum_{d,\alpha,\beta,\gamma}\frac{d}{\delta}
(1+\frac{V}{\beta}+
\frac{V^{2}}{\beta^{2}\gamma})\\
&\ll& U^{1+\ep}\sum_{d,\alpha,\beta,\gamma}\frac{d}{\delta}(V+
\frac{V^{2}}{\beta^{2}\gamma}).
\end{eqnarray*}
We now set
\[f(d)=\prod p^{[e/2]},\]
whence
\[f(d)=\frac{d}{\delta},\]
and
\[g(d)=\prod p^{[e/3]-[(e+1)/2]},\]
for which we claim that
\[g(d)\ge \frac{d}{\delta\beta^{2}\gamma}.\]
To prove the latter it is enough to verify that
\[[\frac{e}{3}]-[\frac{e+1}{2}]\ge e-[\frac{e+1}{2}]
-2[\frac{f+2}{3}]-\max\{e-3[\frac{f+2}{3}]\, ,\, 0\}\]
for $0\le f\le e$, which is an easy exercise.
Moreover, since $\alpha,\beta$ and $\gamma$ all divide $d,$ they take
$O(U^{\ep})$ values 
each.    It follows that 
\[S\ll U^{1+4\ep}\sum_{d\leq U^{2}}\{f(d)V+g(d)V^{2}\}.\]
We now observe that the Dirichlet series $\sum f(d)d^{-\sigma}$ and
$\sum g(d)d^{-\sigma}$ are convergent for $\sigma>1$ and $\sigma>0$ 
respectively, since their Euler products converge.  We may therefore deduce
that
\[\sum_{d\leq U^{2}}f(d)d^{-1-\ep}\ll_{\ep} 1,\]
whence
\[\sum_{d\leq U^{2}}f(d)\leq\sum_{d\leq U^{2}}f(d)(\frac{U^{2}}{d})^{1+\ep}
\ll_{\ep} U^{2+2\ep}.\]
In a similar manner we find that
\[\sum_{d\leq U^{2}}g(d)\ll_{\ep} U^{2\ep}.\]
These bounds suffice for the proof of the lemma, on replacing $\ep$ 
by $\ep/6.$

\section{Theorem 1---Completion of the Proof}

Whenever $w_T(E)\not=0$ we have $\Delta_E\ll T$, and hence $N_E\ll
T$.  It therefore follows from (\ref{2.1}) and (\ref{2.2}) that
\[\frac{1}{{\cal S}(T)}\sum_{E\in{\cal C}}w_{T}(E)r(E)\leq
\frac{\log T}{\log X}+\frac{1}{2}+O(\frac{1+U_1}{\log X})\]
for $X\le T^{2/3}$, where
\begin{eqnarray*}
U_1&=&\frac{1}{{\cal S}(T)}
|\sum_{5\le p\leq X}\frac{\log p}{p}h_{X}(\log p)
\sum_{E\in{\cal C}}w_{T}(E)a_{p}(E)|\\
&=&\frac{1}{{\cal S}(T)}
|\sum_{5\le p\leq X}\frac{\log p}{p}h_{X}(\log p)
\sum_{E\in{\cal C}}w_{T}(E)\sigma_{p}(E)|,
\end{eqnarray*}
by (\ref{3.1a}).
We proceed to show that if $\delta>0$ is fixed, then $U_{1}\ll 1$ when 
$X=T^{2/3-\delta}$.  This suffices for Theorem 1.  

We begin by considering 
the contribution made by the set of singular curves, which 
we denote by ${\cal S}.$  If $E$ is singular then $\Delta_{E}=0,$ whence
\[\sum_{E\in {\cal S}}w_{T}(E)\ll T^{1/6}.\]
Moreover one may verify, using the definition (\ref{3.1}), that 
$|\sigma_p(E)|\leq 1$ for $E\in {\cal S},$ whence
\begin{equation}\label{6.1}
\sum_{P<p\leq 2P}|\sum_{E\in {\cal S}}w_{T}(E)\sigma_p(E)|\ll PT^{1/6}.
\end{equation}

For the non-singular curves $E=E_{r,s}$ we put $r=d^{4}\rho$ and
$s=d^{6}\sigma,$ 
where $d$ is a positive integer and $E_{\rho,\sigma}\in {\cal C}.$  For each
curve $E=E_{r,s}$ we write $E^{*}$ for the corresponding curve
$E_{\rho,\sigma},$ 
so that
\[w_{T}(E)=w_{Td^{-12}}(E^{*}).\]
Moreover $\sigma_p(E)=\sigma_p(E^{*}),$ if $p\nmid d,$ and $\sigma_{p}(E)=0$ 
otherwise.
It follows that
\[\sum_{E\in {\cal D}-{\cal S}}w_{T}(E)\sigma_p(E)=\sum_{d=1}^{\infty}
\sum_{E\in{\cal C}}w_{Td^{-12}}(E)\sigma_p(E) +\theta(T,p),\]
where
\[\theta(T,p)=
-\sum_{d\equiv 0\!\mod{p}}\sum_{E\in{\cal C}}w_{Td^{-12}}(E)\sigma_p(E)
\ll T^{5/6}p^{-10}.\, p^{/2}\]
by the bound (\ref{3.1b}).
The M\"{o}bius inversion formula now yields
\[\sum_{E\in {\cal C}}w_{T}(E)\sigma_p(E)=\sum_{d=1}^{\infty}\mu(d)
\sum_{E\in {\cal D}-{\cal S}}w_{Td^{-12}}(E)\sigma_p(E)-
\sum_{d=1}^{\infty}\mu(d)\theta(Td^{-12},p),\]
and the second sum on the right is
\[\ll \sum_{d=1}^{\infty}(Td^{-12})^{5/6}p^{-19/2}\ll T^{5/6}p^{-19/2}.\]
We therefore see that
\begin{eqnarray}\label{6.2}
\lefteqn{\sum_{5\le p\leq X}\frac{\log p}{p}h_{X}(\log p)
\sum_{E\in{\cal C}}w_{T}(E)\sigma_p(E)}\hspace{1cm}\nonumber\\
&=&\sum_{d=1}^{\infty}\mu(d)
\sum_{5\le p\leq X}\frac{\log p}{p}h_{X}(\log p)
\sum_{E\in{\cal D}-{\cal S}}w_{Td^{-12}}(E)\sigma_p(E)\nonumber\\
&&\hspace{3cm}\mbox{}+
O(\sum_{5\le p\leq X}\frac{\log p}{p}h_{X}(\log p)T^{5/6}p^{-19/2})\nonumber\\
&=&\sum_{d=1}^{\infty}\mu(d)
\sum_{5\le p\leq X}\frac{\log p}{p}h_{X}(\log p)
\sum_{E\in{\cal D}-{\cal S}}w_{Td^{-12}}(E)\sigma_p(E)\nonumber\\
&&\hspace{4cm}\mbox{}+O(T^{5/6}).
\end{eqnarray}
Since ${\cal S}(T)\gg T^{5/6}$ the error term is satisfactory for the
desired bound $U_1\ll 1$.

If $d\gg T^{1/12}$ then $w_{Td^{-12}}(E)$ will vanish.  
Thus we may restrict the sum over $d$ to
the interval $d\ll T^{1/12}$.  We split this range at $d=d_0$, with a
value of $d_0$ to be specified in due course, see (\ref{6.4}).

When $d\ge d_0$ we write
\begin{eqnarray}\label{6.3}
\lefteqn{\sum_{5\le p\leq X}\frac{\log p}{p}h_{X}(\log p)
\sum_{E\in{\cal D}-{\cal S}}w_{Td^{-12}}(E)\sigma_p(E)}\hspace{1cm}\nonumber\\
&=&
\sum_{E\in{\cal D}-{\cal S}}w_{Td^{-12}}(E)
\sum_{5\le p\leq X}\frac{\log p}{p}h_{X}(\log p)\sigma_p(E).
\end{eqnarray}
To estimate the inner sum we put $E=E_{r,s}$ and we let
$r= f^4\rho$ and $s= f^6\sigma$ with $E_{\rho,\sigma}\in 
{\cal C}$. For convenience we set $E^*=E_{\rho,\sigma}$ as before.
Thus $\sigma_p(E)=0$ for $p| f$ and
\[\sigma_p(E)=\sigma_p(E^*)=a_p(E^*)\]
otherwise.  It follows that
\begin{eqnarray*}
\sum_{5\le p\leq X}\frac{\log p}{p}h_{X}(\log p)\sigma_p(E)&=&
\sum_{5\le p\leq X}\frac{\log p}{p}h_{X}(\log p)a_p(E^*)\\
&&\hspace{1cm}\mbox{}-
\twosum{5\le p\leq X}{p| f}\frac{\log p}{p}h_{X}(\log
p)a_p(E^*).
\end{eqnarray*}
The first sum on the right is $O(\log T)^2$ by Lemma 1, while the second
is trivially $O( f)$, since $|a_p(E^*)|\le 2\sqrt{p}$.  Thus the
first sum contributes
\[\ll \sum_{E\in{\cal D}-{\cal S}}w_{Td^{-12}}(E)(\log T)^2
\ll  T^{5/6}d^{-10}(\log T)^2\]
to (\ref{6.3}), and the second
\begin{eqnarray*}
&\ll&\sum_{r}\sum_{s}w_{1}(\frac{r}{(Td^{-12})^{1/3}})
w_{2}(\frac{s}{(Td^{-12})^{1/2}})
\sum_{ f^4|r,\, f^6|s} f \\
&=&\sum_{ f=1}^{\infty} f\sum_{r\equiv 0\!\mod{ f^{4}}}
\sum_{s\equiv 0\!\mod{ f^{6}}}w_{1}(\frac{r}{(Td^{-12})^{1/3}})
w_{2}(\frac{s}{(Td^{-12})^{1/2}})\\
&\ll&\sum_{ f=1}^{\infty} f\,\frac{T^{5/6}d^{-10}}{ f^{10}}\\
&\ll& T^{5/6}d^{-10}.
\end{eqnarray*}
We therefore see that
\[\sum_{5\le p\leq X}\frac{\log p}{p}h_{X}(\log p)
\sum_{E\in{\cal D}-{\cal S}}w_{Td^{-12}}(E)\sigma_p(E)\ll 
T^{5/6}d^{-10}(\log T)^2.\]
Thus terms with $d\ge d_0$ contribute
$O(T^{5/6}d_0^{-9}(\log T)^2)$ to (\ref{6.2}).  On choosing
\begin{equation}\label{6.4}
d_0=\log T,
\end{equation}
say, we see that this is $O(T^{5/6})$, which is satisfactory.

For the values $d<d_0$ we observe that
\begin{eqnarray}\label{6.5}
\lefteqn{\sum_{P<p\leq 2P}
|\sum_{E\in {\cal D}-{\cal S}}w_{Td^{-12}}(E)\sigma_p(E)|}\hspace{1cm}
\nonumber\\
&\leq &\sum_{P<p\leq 2P}|\sum_{E\in {\cal D}}w_{Td^{-12}}(E)\sigma_p(E)|
+O(P(Td^{-12})^{1/6}),
\end{eqnarray}
by a second application of (\ref{6.1}).
According to Lemma 2 the inner sum is
\[\ll P^{\ep}(P^{1/2}T^{5/6}d^{-10}+P^{3/2}T^{1/2}d^{-6}+P^{2}T^{1/6}d^{-2}+
+P^{7/2}T^{-1}d^{12}),\]
whence
\begin{eqnarray*}
\lefteqn{\sum_{P<p\leq 2P}
|\sum_{E\in {\cal D}-{\cal S}}w_{Td^{-12}}(E)\sigma_p(E)|}\\
&\ll&
P^{\ep}(P^{1/2}T^{5/6}d^{-10}+P^{3/2}T^{1/2}d^{-6}+P^{2}T^{1/6}d^{-2}+
P^{7/2}T^{-1}d^{12}),
\end{eqnarray*}
since the error term $O(PT^{1/6}d^{-2})$ in (\ref{6.5}) is majorized by
the term 
\[P^{\ep}.\,P^2T^{1/6}d^{-2}\]
above.  It follows that
\begin{eqnarray*}
\lefteqn{\sum_{5\le p\leq X}\frac{\log p}{p}h_{X}(\log p)
|\sum_{E\in {\cal D}-{\cal S}}w_{Td^{-12}}(E)\sigma_p(E)|}\\
&\ll&
T^{5/6}d^{-10}+X^{2\ep}(X^{1/2}T^{1/2}d^{-6}+XT^{1/6}d^{-2}+
X^{5/2}T^{-1}d^{12}),
\end{eqnarray*}
whence the terms with $d<d_0$ contribute
\[\ll T^{5/6}+X^{2\ep}T^{\ep}(X^{1/2}T^{1/2}+XT^{1/6}+X^{5/2}T^{-1})\]
to (\ref{6.2}).  If we choose $X=T^{2/3-\delta}$, and take $\ep$
sufficiently small in terms of $\delta$, all these terms will be
$O(T^{5/6})$.  This is also satisfactory for the desired bound $U_1\ll
1$.  The proof of Theorem 1 is therefore complete.

\section{Proof of Theorem 2}
To establish Theorem 2 we combine (\ref{2.1}) with the estimate 
(\ref{2.3}) to show that
\[r(E)\leq \frac{\log N_{E}}{\log X}+
\frac{2}{\log X}U_{1}(E,X)+2 +O(\frac{1}{\log X}).\]
It will be convenient to remove the first few primes from the sum
$U_1(E,X)$, so we shall write
\[U_{1}(E,X)=U(E,X)+O(1),\]
where
\[U(E,X)=\sum_{100<p\leq X}c_{p}(E)h_{X}(\log p)\log p,\]
say.  Then, since the curves under consideration have $N_{E}\ll T$, 
we deduce that
\[r(E)\leq 2+\frac{\log T}{\log X}+\frac{2}{\log X}U(E,X)
+O(\frac{1}{\log X}).\]
Consequently, if $X\geq X_{0},$ where $X_{0}$ is a sufficiently large 
absolute constant, and if 
\begin{equation}\label{7.1}
r(E)\geq R\geq 3+2\frac{\log T}{\log X},
\end{equation}
then
\[|U(E,X)|\geq \frac{1}{2}\log T.\]

We complete the proof of the theorem by 
estimating moments of the sum $U(E,X).$  Under the hypothesis (\ref{7.1})
we see that
\[\#\{E\in{\cal C}(T): r(E)\geq R\}(\frac{1}{2}\log T)^{2k}\leq 
\sum_{E\in{\cal C}(T)}|U(E,X)|^{2k},\]
for any positive integer $k.$  We now set
\[V(E,X)=\sum_{100<p\leq X}\frac{\log p}{p}h_{X}(\log p)\sigma_{p}(E)\]
for any $E\in{\cal D}(T)$, so that $U(E,X)=V(E,X)$ whenever
$E\in{\cal C}(T)$, by (\ref{3.1a}).  We then have
\begin{equation}\label{7.2}
\#\{E\in{\cal C}(T): r(E)\geq R\}(\frac{1}{2}\log T)^{2k}\leq 
\sum_{E\in{\cal D}(T)}|V(E,X)|^{2k},
\end{equation}
for any positive integer $k.$  We note that $V(E,X)$ is in fact real,
and expand $|V(E,X)|^{2k}$ by the 
multinomial theorem.  This gives
\begin{equation}\label{7.3}
\sum_{E\in{\cal D}(T)}\sum_{{\bf e}}C({\bf e})F({\bf e}),
\end{equation}
where
\[F({\bf e})=
\prod_{100<p\leq X}\{{\log p}{p}h_{X}(\log p)\sigma_p(E)\}^{e_{p}}.\]
Here ${\bf e}$ runs over vectors with one non-negative integer component
$e_{p}$ for each prime $p\in (100,X],$ and such that $\sum e_{p}=2k.$  Moreover
the multinomial coefficients $C({\bf e})$ are given by
\[C({\bf e})=\frac{(2k)!}{\prod e_{p}!}.\]
We divide the terms in (\ref{7.3}) 
into two classes.  Type I terms will be those
for which every exponent $e_{p}$ satisfies either $e_{p}=0$ or $e_{p}\geq 2.$
The remaining terms will be type II terms.

We begin by considering type I terms.  Since 
$|\sigma_{p}(E)|\leq 2p^{-1/2}$ by (\ref{3.1b}), we
have
\[|F({\bf e})|\leq \prod_{p}(\frac{2h_{X}(\log p)\log p}{\sqrt{p}})^{f_{p}},\]
where 
\[f_{p}=\left\{\begin{array}{cc} 0, & e_{p}=0,\\ 2, & e_{p}\geq 2.
\end{array}\right.\]
Here we use the fact that
\[\frac{2h_{X}(\log p)\log p}{\sqrt{p}}\leq 1\]
for $p>100.$
Moreover $C({\bf e})\leq (2k)!$ for every vector ${\bf e}.$ Thus the terms for
which exactly $j$ primes have $f_{p}=2$ can contribute at most
\[ \frac{(2k)!}{j!}S^{j}\]
to (\ref{7.2}), where
\[S=\sum_{100<p\leq X}\frac{(2h_{X}(\log p)\log p)^{2}}{p}.\]
We now observe that, with our choice of $h_{X},$ we have
\[S=\frac{\log^{2}X}{3}+O(\log X).\]
Thus, if $k\leq \log X$ with $X$ sufficiently large, the contribution 
to (\ref{7.2}) from all type I terms is at most
\begin{equation}\label{7.4}
\ll T^{5/6}\sum_{j\leq k}\frac{(2k)!}{j!}S^{j}\ll
T^{5/6}\frac{(2k)!}{k!}(\frac{\log^{2}X}{3})^{k}\ll 
T^{5/6}(\frac{4k\log^{2}X}{3e})^{k}.
\end{equation}

We turn now to the type II terms.  We begin by recalling the definition
\[\sigma_p(E_{r,s})=
-\tau_{p}^{-1}\sum_{t,x\!\mod{p}}(\frac{t}{p})e_{p}(tx^{3}+txr+ts).\]
When we sum over $E\in{\cal D}(T)$ we have therefore to estimate
\[\sum_{E_{r,s}\in{\cal D}(T)}\hspace{1mm}\sum_{t_{1},\ldots,t_{2k}}
\hspace{1mm}\sum_{x_{1},\ldots,x_{2k}}\hspace{1mm}
\prod_{i=1}^{2k}\{(\frac{t_{i}}{p_{i}})
e_{p_{i}}(t_{i}x_{i}^{3}+t_{i}x_{i}r+t_{i}s)\},\]
where $t_{i}$ and $x_{i}$ run modulo $p_{i},$ and the primes 
$p_{1},\ldots,p_{2k}$ include at least one value, $p^{*}$ say, which is 
not repeated.  We bound the above expression as
\[\ll (\prod p_{i})\sum_r\sum_{t_{1},\ldots,t_{2k}}|\sum_{s}
e(s\{\frac{t_{1}}{p_{1}}+\ldots+\frac{t_{2k}}{p_{2k}}\})|,\]
in which $t_{i}$ runs over $1,\ldots,p_{i}-1$, and the function $e(x)$
is given by $e(x)=\exp(2\pi ix)$.  The summation conditions on $r$ and
$s$ are given by
\[|r|\le T^{1/3},\;\;\; |s|\le T^{1/2},\;\;\mbox{and}\;\;
4r^3+27s^2\not=0.\]
We proceed to examine the innermost sum.
We write
\[\frac{t_{1}}{p_{1}}+\ldots+\frac{t_{2k}}{p_{2k}}=\alpha\]
and note that $\alpha$ cannot be an integer, since its denominator must be 
divisible by $p^{*}.$  It follows that
\[||\alpha||\geq (\prod p_{i})^{-1},\]
whence
\[|\sum_{s}e(s\alpha)|\leq 2+\frac{1}{|\sin\pi\alpha|}\ll \frac{1}{||\alpha||}
\leq \prod p_{i}.\]
Here we have allowed for the fact that, for a given value of $r$, the
variable $s$ runs over all integers in the interval
$[-T^{1/2},T^{1/2}]$ with at most 2 exceptions.  Since $r$ takes
$O(T^{1/3})$ values
we therefore see that 
\[\sum_{E\in{\cal D}(T)}F({\bf e})\ll 
T^{1/3}\prod_{100<p\leq X}\{\frac{p^2\log p}{p|\tau_p|}\}^{e_{p}}\]
for each type II term, whence the total contribution to (\ref{7.2}) is
\begin{equation}\label{7.5}
\ll T^{1/3}(\sum_{100<p\leq X}p^{1/2}\log p)^{2k}\ll T^{1/3}X^{3k}.
\end{equation}

In view of (\ref{7.2}) and the estimates (\ref{7.4}) and (\ref{7.5}) 
we find that
\[\#\{E\in{\cal C}(T): r(E)\geq R\}(\frac{1}{2}\log T)^{2k}\ll 
T^{5/6}(\frac{4k\log^{2}X}{3e})^{k}+T^{1/3}X^{3k},\]
for $X\geq X_{0},$ subject to the conditions
\[R\geq 3+2\frac{\log T}{\log X}\]
and $k\leq \log X.$  Note here that $X_{0}$ is independent of $k.$
We therefore choose 
\[X=T^{1/6k},\]
whence
\[\#\{E\in{\cal C}(T): r(E)\geq R\}\ll (27ek/4)^{-k}T^{5/6},\]
for $R\geq 3+12k$ and $T^{1/6k}\geq X_{0}.$   We take
$k=[\frac{R-3}{12}]$ and write $j=R/12$, so that $k\le j\le k+O(1)$.
Then for any positive constants $a>b$ we will have $(ak)^{-k}\ll
(bj)^{-j}$ if $k$ is large enough.  Since $27e/4>18$ we conclude that
\begin{equation}\label{7.6}
\#\{E\in{\cal C}(T): r(E)\geq R\}\ll (3R/2)^{-R/12}T^{5/6},
\end{equation}
if $R$ is large enough and $R\ll \log T.$  However (\ref{7.6}) 
is trivially true
for bounded values of $R.$ Moreover for 
\[R=[11\frac{\log T}{\log\log T}]\]
we may already conclude from (\ref{7.6}) that
\[\#\{E\in{\cal C}(T): r(E)\geq R\}=o(1),\]
so that there can be no curves with
\[r(E)\geq 11\frac{\log T}{\log\log T}\]
for large enough $T.$
This completes the proof of Theorem 2.

\section{Theorem 3---Preliminary Sieving}

The condition $D\in{\cal T}^{\pm}$ is distinctly awkward to work with, and
our first task is therefore to replace it with something more manageable.
When $N$ is odd we begin by decomposing ${\cal T}^{\pm}$ according to the 
power of $2$ dividing $D.$  Of course, if $N$ is even then $D$ is 
automatically odd.  We now write $D=\delta 2^{e}n$ with $\delta=+1$ or
$-1$ and $n$ odd, and we decompose 
${\cal T}^{\pm}$ further according to the residue class of $n$ modulo $8.$
This produces a collection of triples $(k,\delta,e),$ in which 
\[k=1,3,5\;\mbox{or}\; 7,\;\;\;\delta=+1\;\mbox{or}\; -1,\;\;\;
\mbox{and}\;\;\; e=0,2,\;\mbox{or}\; 3,\]
and such that 
${\cal T}^{\pm}$ is a disjoint union of certain of the sets
\[ \{D=\delta 2^{e}n: w_D=\pm 1,\,\mu^2(n)=1, n\equiv k\!\mod{8}\}.\]
We shall prove the analogue of Theorem 3 for these sets,
assuming that the weight function $w$ is supported
on a compact subset of $(-\infty,0)$ for $\delta=-1,$ and $(0,\infty)$ for
$\delta=+1.$  Theorem 3 itself will then follow.  Henceforth we shall regard
the triple $(k,\delta,e)$ and the sign $\pm1$ as fixed, and 
for any positive odd integer $n$ we shall write $D=D(n)=
\delta 2^{e}\hat{n}$, where $\hat{n}$ is the square-free kernel of $n$.
We also set
\[W(n/T)=w(\delta 2^{e}n/T),\]
\[{\cal F}=\{n\in\N: w_D=\pm 1,\,(n,N)=1,\,\mu^2(n)=1,\, 
n\equiv k\!\mod{8}\},\]
and
\[{\cal R}(T)=\sum_{n\in{\cal F}}W(n/T).\]
We have therefore to show that
\[\frac{1}{{\cal R}(T)}\sum_{n\in{\cal F}}W(n/T)r(E_{D})\leq 
\frac{3}{2}+o(1),\]
as $T\rightarrow\infty.$

We turn now to the condition that $n$ must be square-free.  We define
\[P=\prod_{\begin{array}{c} {\scriptstyle 2<p\leq \log\log T}\\
{\scriptstyle p\nmid N} \end{array}}p,\]
and we set
\[ X(n)=\sum_{d|P,\,d^{2}|n}\mu(d),\]
so that $X(n)=0$ if $n$ is divisible by the square of a prime 
$p\leq\log\log T,$ and $X(n)=1$ otherwise.  It follows that  
\[\sum_{n\in{\cal F}}W(n/T)r(E_{D})\leq 
\sum_{n\in{\cal G}}X(n)W(n/T)r(E_{D}),\]
where
\[{\cal G}=\{n\in\N: w_D=\pm 1,\,(n,N)=1,\, n\equiv k\!\mod{8}\}.\]
Moreover it is a straightforward matter to demonstrate the asymptotic formula
\[\sum_{n\in{\cal G}}X(n)W(n/T)\sim {\cal R}(T),\]
since ${\cal R}(T)\gg T$ if ${\cal F}$ is non-empty.
It therefore suffices to establish the estimate
\begin{equation}\label{8.1}
\sum_{n\in{\cal G}}X(n)W(n/T)r(E_{D})\leq (\frac{3}{2}+o(1))
\sum_{n\in{\cal G}}X(n)W(n/T).
\end{equation}

The proof of Theorem 3 now hinges on the following lemma.
\begin{lemma}
Let $2\leq X\leq T^{2-\ep},$ where $\ep$ is a positive constant.  Then, 
assuming the 
Riemann Hypothesis for all the L-functions $L_{D}(s),$ we have
\[\sum_{n\in{\cal G}}W(n/T)U_{1}(E_{D(n)},X)\ll T.\]
\end{lemma}

We conclude this section by demonstrating how (\ref{8.1}) may be deduced from
Lemma 4.  We have
\begin{eqnarray}\label{8.2}
\lefteqn{\sum_{n\in{\cal G}}X(n)W(n/T)U_{1}(E_D,X)}\hspace{2cm}\nonumber\\
&=&
\sum_{d|P}\mu(d)\sum_{n\in{\cal G}: d^{2}|n}W(n/T)U_{1}(E_{D(n)},X)
\nonumber\\
&=&\sum_{d|P}\mu(d)\sum_{m\in{\cal G}}W(d^{2}m/T)U_{1}(E_{D(m)},X),
\end{eqnarray}
on replacing $n$ by $d^{2}m.$  (Notice that $D(d^2 m)=D(m)$, 
and that $m\in{\cal G}$ if and only if $d^{2}m\in{\cal G}.$)
However
\[d\leq P=\exp\{O(\log\log T)\}\ll T^{\ep}.\]
Thus if $X\leq T^{2-5\ep}$ we have
\[X\leq (\frac{T}{d^{2}})^{2-\ep},\]
for any $\ep>0.$  We may therefore apply 
Lemma 4 to the inner sum in (\ref{8.2}), giving
\[\sum_{n\in{\cal G}}X(n)W(n/T)U_{1}(E_D,X)\ll 
\sum_{d|P}\frac{T}{d^{2}}\ll T\]
for $X= T^{2-5\ep}$.  We now feed (\ref{2.1}) and (\ref{2.2a}) into
the left-hand side of (\ref{8.1})
and observe that 
\[N_{E_D}\ll D^2\ll T^2\]
since the curve $E$ is fixed.  This produces the required bound (\ref{8.1}).
Thus to complete the proof of Theorem 3 it will suffice to establish
Lemma 4.

\section{Further Simplifications}
In this section we shall simplify the expression occurring in Lemma 4.
We begin by noting that
\[a_p(E_D)=(\frac{D}{p})a_p(E)\]
for primes $p\nmid ND$, and  we therefore classify the odd 
primes $p$ according to their residue modulo $8,$ which enables us to write
\[(\frac{D}{p})=\eta_{h}(\frac{n}{p}),\;\;\mbox{for}\;\; p\equiv h\!\mod{8},\; 
p\nmid n,\]
where $\eta_{h}$ may depend on $\ep,\delta$ and $e$ as well as on $h.$ 
Thus
\[U_{1}(E_D,X)=-\sum_{h}\eta_{h}\sum_{p\equiv h\!\mod{8}}\beta_{p}(\frac{n}{p})
+O(\sum_{p|Nn}\frac{\log p}{\sqrt{p}}),\]
where we have introduced the shorthand
\[\beta_{p}=\left\{\begin{array}{cc} \frac{\log p}{p}h_{X}(\log p)a_{p}(E), &
p\ge 5,\\ 0, & p=2,\,3.\end{array}\right.\]
Since $E$ is fixed we have
\[\sum_{p|N}\frac{\log p}{\sqrt{p}}\ll 1,\]
so that these terms contribute $O(T)$ in Lemma 4.  Moreover
\begin{eqnarray*}
\sum_{n}W(n/T)\sum_{p|n}\frac{\log p}{\sqrt{p}} &=&
\sum_{p}\frac{\log p}{\sqrt{p}}\sum_{p|n}W(n/T)\\
&\ll&\sum_{p}\frac{\log p}{\sqrt{p}}\frac{T}{p}\\
&\ll& T,
\end{eqnarray*}
which is also satisfactory.  The condition $p\equiv h\mod{8}$ may 
be picked out by using an appropriate combination of the characters
\[(\frac{a}{p}),\;\;\;a=1,-1,2,-2.\]
For the proof of Lemma 4 it therefore suffices
to show that
\[\sum_{n\in{\cal G}}W(n/T)U(an)\ll T,\]
for $a=1,-1,2,-2,$ where
\[U(m)=\sum_{p}\beta_{p}(\frac{m}{p}).\]

We turn now to the condition $n\in{\cal G}.$  Since $(n,N)=1,$ the
root number
$w_D$ differs from $(\frac{N}{n})$ by a factor depending on 
$N,k,\delta$ and $e$ only, in view of (\ref{1.1}).  
We can therefore pick out the conditions 
$w_D=\pm 1$ and $n\equiv k\mod{8}$ by introducing a suitable 
combination of factors $(\frac{N}{n})$, $(\frac{2}{n})$ and 
$(\frac{-1}{n}).$  We deduce that
it is sufficient, for the proof of Lemma 4, to establish the estimate
\[\sum_{(n,2N)=1}\psi(n)W(n/T)U(an)\ll T,\]
where $\psi(n)$ is a real primitive character of conductor dividing $8N$.
Since
\begin{eqnarray*}
\sum_{(n,2N)=1}\psi(n)W(n/T)U(an)&=&
\sum_{d|2N}\mu(d)\sum_{d|n}\psi(n)W(n/T)U(an)\\
&=&\sum_{d|2N}\mu(d)\psi(d)\sum_{m=1}^{\infty}\psi(m)W(dm/T)U(adm),
\end{eqnarray*}
we conclude as follows.
\begin{lemma}
In order to establish Lemma 4 it suffices to show that
\[\sum_{n=1}^{\infty}\psi(n)W(n/T)U(rn)\ll_{r,E} T,\]
for $2\leq X\leq T^{2-\ep},$ and for every $r\not=0.$
\end{lemma}

\section{Character Sums}
We now have to examine
\begin{equation}\label{10.1}
\sum_{n=1}^{\infty}\psi(n)W(n/T)U(rn)=\sum_{p}\beta_{p}(\frac{r}{p})
\sum_{n}W(n/T)\psi_{p}(n),
\end{equation}
where
\begin{equation}\label{10.2}
\psi_{p}(n)=\psi(n)(\frac{n}{p}).
\end{equation}
We shall denote the sum on the left of (\ref{10.1}) by $\Sigma.$
The primes for which $p|N,$ contribute a total $O(T)$ to $\Sigma.$  
For the remaining primes $\psi_{p}$ is primitive.   We write $\Sigma_{p}$ for
the inner sum on the right of (\ref{10.1}), and we denote the conductor 
of $\psi_{p}$
by $q.$  Thus $q=bp,$ say, where $b$ is the conductor of $\psi$.  
Moreover $b|8N.$
We proceed to decompose $\Sigma_{p}$ by dividing the values of $n$ into
congruence classes $n\equiv j\mod{q},$ whence
\[\Sigma_{p}=\sum_{j\!\mod{q}}\psi_{p}(j)\sum_{m=-\infty}^{\infty}
W(\frac{j+qm}{T}).\]
On applying the Poisson summation formula we obtain
\begin{eqnarray*}
\Sigma_{p}&=&\sum_{j\!\mod{q}}\psi_{p}(j)\sum_{m=-\infty}^{\infty}
\frac{T}{q}e(\frac{mj}{q})\hat{W}(\frac{Tm}{q})\\
&=&\frac{T}{q}\sum_{m=-\infty}^{\infty}\hat{W}(\frac{Tm}{q})
\sum_{j\!\mod{q}}\psi_{p}(j)e(\frac{mj}{q}),
\end{eqnarray*}
where $e(x)=\exp\{2\pi ix\}$ as usual.  On writing $G(p)$ for the Gauss sum
\[\sum_{j\!\mod{q}}\psi_{p}(j)e(\frac{j}{q}),\]
we have
\[\sum_{j\!\mod{q}}\psi_{p}(j)e(\frac{mj}{q})=G(p)\psi_{p}(m),\]
so that
\[\Sigma_{p}=T\frac{G(p)}{q}\sum_{m=-\infty}^{\infty}
\hat{W}(\frac{Tm}{q})\psi_{p}(m).\]
Since $\psi_{p}(0)=0,$ we therefore conclude that
\begin{equation}\label{10.3}
\Sigma=\frac{T}{b}\sum_{m\not=0}\sum_{p\nmid N}\frac{G(p)}{p}\beta_{p}
\hat{W}(\frac{Tm/b}{p})\psi_{p}(m)(\frac{r}{p})+O(T).
\end{equation}

It is instructive to examine the trivial estimate for $\Sigma$ at this
stage.  Since the function $W$ is supported on a compact subset of
$(0,\infty),$ and is three times differentiable, we have
\begin{equation}\label{10.4}
\hat{W}(x)\ll \min\{1\, ,\,|x|^{-3}\}.
\end{equation}
Thus, on using the bounds $G(p)\ll p^{1/2}$ and $\beta_{p}\ll p^{-1/2}\log p,$
we find that
\begin{eqnarray*}
\Sigma&\ll& T+T\sum_{m\not=0}\;\sum_{p\leq X}\frac{\log p}{p}
\frac{p^{3}}{T^{3}|m|^{3}}\\ 
&\ll & T+T^{-2}X^{3}.
\end{eqnarray*}
This therefore suffices for an analogue of Lemma 5 in 
which $X$ may be as large 
as $T^{1-\ep}.$  One would then obtain a version of Theorem 3 with a constant 
$\frac{5}{2}$ in place of $\frac{3}{2}.$  Such an improvement of Goldfeld's
bound was mentioned by Brumer \cite[p. 445]{brum}, 
although it is clear that the argument 
intended by Brumer was a relatively minor modification of that used by 
Goldfeld.

Our sharper estimate for $\Sigma$ stems from a non-trivial bound for
the inner sum in (\ref{10.3}).  To obtain this we call on the following
`Prime Number Theorem' for twisted curves $E_{D}.$
\begin{lemma}
If $L_{D}(s)$ satisfies the Riemann Hypothesis
we have
\[\sum_{p\leq x}\frac{a_{p}(E)}{p}\chi_{D}(p)\log p\ll x^{\ep}|D|^{\ep}\]
for any $\ep>0,$ where the implied constant depends at most on $E$ and $\ep.$
\end{lemma}

Here the reader should recall that $\chi_D$ is the real primitive
character associated to the quadratic field $\Q(\sqrt{D})$.  (When
$D=1$ we take $\chi_D$ to be the trivial character.) 
We shall prove Lemma 6 in the next section.  Notice that the lemma does not
assume that $D$ and $N$ are coprime.

To apply Lemma 6 to (\ref{10.3}) we observe that
\[G(p)=\psi(p)(\frac{b}{p})\tau(\psi)\tau_p=
C_{b}\psi'(p)\sqrt{p}(1-i(\frac{-1}{p})),\]
by the usual evaluation of Gauss sums.  Here $C_{b}$ is a constant depending on
$b$ only, and $\psi'$ is a real character whose modulus divides $8N.$  In view
of the definition (\ref{10.2}) of $\psi_{p},$ it follows that there is a
real character
$\psi_{1}$ whose modulus divides $8Nmr,$ such that
\begin{eqnarray*}
\lefteqn{\sum_{p\nmid N}\frac{G(p)}{p}\beta_{p}
\hat{W}(\frac{Tm/b}{p})\psi_{p}(m)(\frac{r}{p})}\hspace{2cm}\\
&\ll&
|\sum_{p\nmid 30N}\frac{a_{p}(E)}{p^{3/2}}(\log p)h_{X}(\log p)
\hat{W}(\frac{Tm/b}{p})\psi_{1}(p)|.
\end{eqnarray*}
We now wish to replace $\psi_{1}$ by the primitive character $\chi_{\Delta}$
which induces it.  Here $\Delta$ is a fundamental discriminant and 
$\Delta|8Nmr.$  This process will introduce an error which contributes
\[\ll T\sum_{m\not=0}\;
\sum_{p|30Nmr}\frac{\log p}{p}\frac{p^{3}}{T^{3}|m|^{3}}\]
to $\Sigma,$ by (\ref{10.4}).
The primes dividing $30Nr$ provide at most
\[\ll T\sum_{m\not=0}\;\sum_{p|30Nr}\frac{\log p}{p}\frac{p^{3}}{T^{3}|m|^{3}}
\ll T^{-2},\]
and the primes $p|m$ yield a total
\begin{eqnarray*}
\ll T\sum_{m\not=0}\hspace{2mm}\sum_{p|m,\, p\leq X}
\frac{\log p}{p}\frac{p^{3}}{T^{3}|m|^{3}}
&\ll & T^{-2}\sum_{p\leq X}p^2\log p\sum_{m\not=0,\,p|m}|m|^{-3}\\
&\ll & T^{-2}\sum_{p\leq X}p^2(\log p).p^{-3}\\
&\ll& T^{-2}\log X.
\end{eqnarray*}
Both these contributions are satisfactory, and we conclude 
from (\ref{10.3}) that
\begin{equation}\label{10.5}
\Sigma\ll T+T\sum_{m\not=0}|\sum_{p}\frac{a_{p}(E)}{p^{3/2}}(\log p)
h_{X}(\log p)\hat{W}(\frac{Tm/b}{p})\chi_{\Delta}(p)|.
\end{equation}

We shall bound the sum over $p$ by using partial summation together with
the estimate
\[\sum_{p\leq x}\frac{a_{p}(E)}{p}\chi_{\Delta}(p)\log p\ll x^{\ep}|m|^{\ep}\]
which follows from Lemma 6.  In analogy to (\ref{10.4}) we have
\[\frac{d}{dx}\hat{W}(x)\ll \min\{1\, ,\,|x|^{-3}\}.\]
We deduce that
\begin{eqnarray*}
\lefteqn
{\int_{2}^{X}|\frac{d}{dt}\{t^{-1/2}h_{X}(\log t)
\hat{W}(\frac{Tm/b}{t})\}|dt}\hspace{4cm}\\
&\ll & \int_{2}^{X}t^{-3/2}\min\{1\, ,\, (\frac{T|m|}{t})^{-2}\}dt\\
&\ll & (T|m|)^{-1/2}\min\{1\, ,\, (\frac{X}{T|m|})^{3/2}\}.
\end{eqnarray*}
We therefore conclude, on summing by parts, that
\begin{eqnarray*}
\lefteqn{\sum_{p}\frac{a_{p}(E)}{p^{3/2}}(\log p)h_{X}(\log p)
\hat{W}(\frac{Tm/b}{p})\chi_{\Delta}(p)}\hspace{2cm}\\
&\ll& (X|m|)^{\ep}(T|m|)^{-1/2}\min\{1\, ,\, (\frac{X}{T|m|})^{3/2}\}.
\end{eqnarray*}
In view of (\ref{10.5}) we now have
\begin{eqnarray*}
\Sigma&\ll& T+T^{1/2}X^{\ep}\sum_{m\not=0}|m|^{-1/2+\ep}
\min\{1\, ,\, (\frac{X}{T|m|})^{3/2}\}\\
&\ll & T+X^{1/2+2\ep},
\end{eqnarray*}
from which Lemma 5 follows on redefining $\ep.$

\section{Proof of Lemma 6}

To prove Lemma 6 we shall apply Lemma 1 to the curve $E_{D}$,
taking the function $k$ to be 
\[k(t)=\frac{Xh(t)-(X-1)h(\frac{t}{1-X^{-1}})}{\log^{2}X}.\]
Since 
\[\hat{k}(t)=\frac{X}{\log^{2}X}
\frac{\sin^{2}(\pi t)-\sin^{2}(\pi(1-X^{-1})t)}{\pi^{2}t^{2}}\]
the hypothesis of Lemma 1 is satisfied for any $\delta>0$.  Moreover
one readily finds that $||k||_{\infty}\ll (\log X)^{-2}$ and that
\[\hat{k}(t)\ll\frac{X}{t^2(\log x)^2}\min\{1\,,\,|t|\}
\min\{1\,,\,\frac{|t|}{X}\},\]
whence
\[\left|\left|(1+|t|)^{1+\delta}\hat{k}(t)\right|\right|_{\infty}\ll
\frac{X^{\delta}}{(\log X)^2}.\]
It follows that
\[\sum_{p\leq X}\frac{\log p}{p}k(\frac{\log p}{\log X})a_{p}(E_{D})
\ll X^{\delta}\log |D|.\]
However
\[k(\frac{\log p}{\log X})=\left\{\begin{array}{cc} (\log X)^{-2}, & 
p\leq X^{1-1/X},\\ O((\log X)^{-2}), & X^{1-1/X}\leq p\leq X,
\end{array}\right.\]
whence
\begin{eqnarray*}
\lefteqn{\sum_{p\leq X}
\frac{\log p}{p}k(\frac{\log p}{\log X})a_{p}(E_{D})}\hspace{2cm}\\
&=&(\log X)^{-2}\sum_{p\leq X}\frac{\log p}{p}a_{p}(E_{D})
+O(X^{-1/2}\log^{2}X).
\end{eqnarray*}
It therefore follows, on choosing $\delta=\ep/2$, that
\[\sum_{p\leq X}\frac{\log p}{p}a_{p}(E_{D})\ll
X^{\ep}|D|^{\ep}.\]
To complete the proof of Lemma 6 it remains to observe that the only
primes for which $a_p(E_D)$ can differ from $a_p(E)\chi_D(p)$ are, possibly,
those for which $p|30ND$.  Since $N$ is fixed,
these contribute $O(|D|^{\ep})$, which is satisfactory.

\section{Deduction of Theorem 4}

To prove Theorem 4 we begin by observing that
\[\sum_{D\in{\cal T}^{+},\, R(E_{D})=0}w(D/T)\geq
\sum_{D\in{\cal T}^{+},\, r(E_{D})=0}w(D/T),\]
by (\ref{1.2}).  Moreover, Theorem 3 yields
\begin{eqnarray*}
\sum_{D\in{\cal T}^{+},\, r(E_{D})=0}w(D/T)&=&
\sum_{D\in{\cal T}^{+}}w(D/T)-\sum_{D\in{\cal T}^{+},\, r(E_{D})\ge 2}w(D/T)\\
&\ge & {\cal W}^{+}(T)-
\sum_{D\in{\cal T}^{+},\, r(E_{D})\ge 2}w(D/T)\frac{r(E_D)}{2}\\
&=&{\cal W}^{+}(T)-
\sum_{D\in{\cal T}^{+}}w(D/T)\frac{r(E_D)}{2}\\
&\ge &{\cal W}^{+}(T)-\frac{1}{2}\{\frac{3}{2}+o(1)\}{\cal W}^{+}(T)\\
&=&\{\frac{1}{4}+o(1)\}{\cal W}^{+}(T),
\end{eqnarray*}
since $r(E_D)/2\ge 1$ whenever $r(E_D)\ge 2$.

Similarly we have
\[\sum_{D\in{\cal T}^{-},\, R(E_{D})=1}w(D/T)\geq
\sum_{D\in{\cal T}^{-},\, r(E_{D})=1}w(D/T),\]
and
\begin{eqnarray*}
\sum_{D\in{\cal T}^{-},\, r(E_{D})=1}w(D/T)&=&
\sum_{D\in{\cal T}^{-}}w(D/T)-\sum_{D\in{\cal T}^{-},\, r(E_{D})\ge 3}w(D/T)\\
&\ge & {\cal W}^{-}(T)-
\sum_{D\in{\cal T}^{-},\, r(E_{D})\ge 3}w(D/T)\frac{r(E_D)-1}{2}\\
&=&{\cal W}^{-}(T)-
\sum_{D\in{\cal T}^{-}}w(D/T)\frac{r(E_D)-1}{2}\\
&\ge &\frac{3}{2}{\cal W}^{-}(T)-
\frac{1}{2}\{\frac{3}{2}+o(1)\}{\cal W}^{-}(T)\\
&=&\{\frac{3}{4}+o(1)\}{\cal W}^{-}(T),
\end{eqnarray*}
as required for Theorem 4.

Mathematical Institute,

24-29, St. Giles',

Oxford

OX1 3LB

England
\bigskip

rhb@maths.ox.ac.uk

\end{document}